\documentclass[12pt]{article}
\usepackage{epsfig, amsfonts, amsmath, amsthm,amsbsy}
\usepackage{showlabels}
\usepackage{color}
\usepackage{subfigure}

\newtheorem{theorem}{Theorem}
\newtheorem{lemma}{Lemma}

\newtheorem{remark}{Remark}
\def\ve{\varepsilon}

\begin{document}

\title{ A numerical method for singularly perturbed convection-diffusion  problems posed on smooth domains}
\author{A. F. Hegarty 
\thanks{Department of Mathematics and Statistics, University of Limerick, Ireland.} \and E.\ O'Riordan \thanks{ School of Mathematical Sciences, Dublin City University, Dublin 9, Ireland.} }

\maketitle

\begin{abstract}

 A finite difference method is constructed  to solve  singularly perturbed convection-diffusion problems posed on smooth domains.  
Constraints are imposed on the data so that only regular exponential boundary layers appear in the solution. 
A  domain decomposition method is used, which  uses a rectangular grid outside the boundary layer and a Shishkin mesh, aligned to the curvature of the outflow boundary, near the  boundary layer.  Numerical results are presented to demonstrate the effectiveness of the proposed numerical algorithm. 
\bigskip


\end{abstract}

\section{Introduction}

The numerical solution of singularly perturbed elliptic problems, of convection-diffusion type, posed on smooth domains presents several challenges. In order to generate a pointwise accurate global approximation to the solution using  piecewise polynomial basis functions,  the grid needs to insert mesh points into the layer regions, where the derivatives of the solution depend inversely on the magnitude of the singular perturbation parameter. To avoid the dimension of the discrete problem depending on the inverse of the singular perturbation parameter, a quasi-uniform discretization of the continuous domain will not suffice \cite{fhmos}.  Outside the layer regions one only needs a coarse mesh and within the layers one requires a fine mesh. However, spurious oscillations will appear on a coarse mesh, unless some particular discretization is used to preserve the inverse monotoniticity of the differential operator \cite{gabriel, bosco}. Finite element discretizations using triangles  are well suited to discretizing domains with  smooth geometries; but  it is difficult to generate an inverse-monotone system matrix  to a singularly perturbed convection-diffusion problem using triangular elements \cite{john1, john2}. Upwinded finite difference operators (or finite volumes \cite{rst2}) can be used to guarantee stability, but it is not easy to work with these constructions  over anisotropic meshes posed on a smooth domain.

 In addition to these complications, our objective  is to design a parameter-uniform numerical method \cite{fhmos}, which will be accurate both for the classical case (where the singular perturbation parameter is large) and the singularly perturbed case (where the singular perturbation parameter is very small), and will also deal with all  the intermediate values of the singular perturbation parameter.  Moreover, given that boundary layers will be present, we are only interested in global approximations that are pointwise accurate at all points in the domain. Hence, although we use a finite difference formulation, our focus is not primarily on the  nodal accuracy of the numerical method, but on the global accuracy of the interpolated approximation generated by the numerical algorithm.

 To circumvent  the difficulties mentioned above,  we employ a domain decomposition algorithm.  The domain is first covered  with a rectangle to generate an initial approximation to the solution. On this rectangle, we simply use an upwinded finite difference operator on a tensor product of uniform meshes. 
This classical method will produce an accurate and stable approximation outside of the boundary layers.  A correction within the layers is generated across a subdomain that is aligned to the curved outflow
boundary of the original domain \cite{christos, christos2}. Across this subdomain, a coordinate transformation is used  and a piecewise-uniform Shishkin  mesh \cite{fhmos}  is employed  in the normal direction to the boundary. 

On any closed domain, there will be characteristic points on the boundary, where the tangent to the boundary is parallel to the convective direction $\vec a$, which is associated with the characteristics of the reduced first order problem. In order to establish the  theoretical error bound in Theorem \ref{Main-Result} below, we impose constraints on the data via three assumptions. The first assumption (\ref{assume0}) prevents characteristic boundary layers forming. The second assumption (\ref{assum1}) ensures that the solution is sufficiently regular for the numerical analysis within the paper to apply. The final assumption (\ref{assum2}) prevents internal layers emerging within the solution.
 The numerical results in \S 4 suggest that these theoretical constraints are excessive, as the numerical method continues to display first order parameter-uniform under significantly weaker constraints. The identification of necessary  data constraints  to retain the error bound in Theorem \ref{Main-Result} remains an open question. 

In \S 2, the continuous problem is discussed and the solution is decomposed into a regular and a singular component. Pointwise bounds on the  derivatives (up to third order)  of these components are deduced. In \S 3, a numerical method is constructed and an asymptotic error bound is deduced in Theorem \ref{Main-Result}. Numerical results for three sample  problems are presented in the final section.  
\vskip0.5cm 

\noindent{\bf Notation}
If $D$ is the domain of some function $f$ and $D \subset D^*$, then throughout we denote the extension of the function to the larger domain by $f^*$. 
In addition,   $f(x,y) = \tilde f (r,t)$, where $(r,t)$ is a co-ordinate system aligned to the boundary of the domain.
 Here and thoughout $C$
denotes a generic constant that is independent of both  the singular perturbation parameter $\ve$ and the discretization parameter $N$.

\section{Continuous problem}

Let $\Omega $ be a two dimensional  domain
with a smooth closed boundary $\partial \Omega$. The origin is located within the domain. As in \cite{kopteva} we introduce 
a local curvilinear coordinate system associated with the boundary. Let the boundary $\partial \Omega$ be parameterized  by
\[
\partial \Omega := \{ (x, y))\vert x =\phi (t), y =\psi (t),\ 0 \leq t \leq T \}, \quad \gamma (t) :=(\phi (t), \psi (t));
\]
where $(\phi (0), \psi(0)) = (\phi (T), \psi(T))$.
As the variable $t$ increases, the boundary points move in an anti-clockwise direction.  
At any point on the boundary, the magnitude of the tangent vector $(\phi ', \psi ')$ is denoted by $\tau (t) $ and the curvature of the boundary by $\vert \kappa (t) \vert $, which are given by 
\[
\tau (t) :=  \sqrt {(\phi' )^2 +  (\psi' )^2} \quad \hbox{and} \quad \kappa (t) :=
\frac{\phi' \psi '' - \psi' \phi ''}{\tau ^3}.
\]
A curvilinear local coordinate system $(r,t)$ is defined by
\begin{equation}\label{Lame}
x=\phi (t) +rn_1(t), \ n_1 := \frac{-\psi '}{\tau}; \qquad y=\psi (t) +rn_2(t), \ n_2 := \frac{\phi '}{\tau}.
\end{equation}
Note that $\vec n = (n_1,n_2)$ is the inward unit normal and 
\[
n_1'=-\phi ' \kappa =-\kappa \tau n_2,\quad  n_2'=-\psi ' \kappa =\kappa \tau n_1.
\] These coordinates are orthogonal in the sense that
\begin{eqnarray*}
(\frac{\partial x }{\partial r} , \frac{\partial y }{\partial r}) \cdot (\frac{\partial x }{\partial t} , \frac{\partial y }{\partial t}) = (n_1,n_2) \cdot (\phi '+rn_1', \psi '+rn_2') =0.
\end{eqnarray*}
In these coordinates 
, the transformed Laplacian will contain no mixed second order derivative and   we have \cite[Lemma 2.1]{kopteva}
\begin{equation}\label{Transformed}
\triangle u = \eta ^{-1} \frac{\partial }{\partial r} \bigl( \eta \frac{\partial \tilde u }{\partial r} \bigr) + 
\zeta  \frac{\partial }{\partial t} \bigl( \zeta \frac{\partial \tilde u }{\partial t} \bigr), \quad \hbox{where} \quad \zeta := \frac{1}{\tau \eta}, 
\eta := 1- \kappa r.
\end{equation}

Consider the singularly perturbed convection-diffusion elliptic problem
 {\footnote{ For nonnegative integers $k$ and all $ v \in C^k(D), \ D \subset \mathbb{R}^2$, we define 
\[
\vert v \vert _{k,D} := \sum _{i+j=k} \sup _{(x,y) \in D} \Bigl
\vert \frac{\partial ^{k} v}{\partial x ^i
\partial y^j} \Bigr \vert, \quad
\Vert v \Vert _{k,D} := \sum _{0\leq j \leq k} \vert v \vert
_{j,D}.
\]
If $D =\bar{\Omega}$ we omit the subscript $D$ and if $k = 0$ we omit the subscript $k$. The space $C^{\gamma}(D ) $ is the set of all functions that are H\"{o}lder continuous of degree $\gamma $ with respect to the
Euclidean norm $\Vert \cdot \Vert _e$.  A function  $f \in C^{\gamma}(D ) $ if
\[
\lceil f \rceil _{0,\gamma , D} = \sup _{\bf{u}
 \neq \bf{v}, \ \bf{u}, \bf{v} \in D} \frac{\vert f({\bf u}) - f({\bf v}) \vert}{\Vert \bf{u}- \bf{v} \Vert _e^\gamma}
\]  is finite.
 The space $C^{k, \varsigma}(D ) $ is
the set of all functions in $C^{k}(D )$ whose derivatives of order
$k$ are H\"{o}lder continuous of degree $\varsigma $.
Also we define
\[
\lceil v \rceil _{k,\varsigma,D} := \sum _{i+j=k} \Bigl \lceil
\frac{\partial ^k v} {\partial x^i \partial y^j} \Bigr \rceil
_{0,\varsigma,D}, \quad \Vert v \Vert _{k,\varsigma,D} := \sum _{0\leq n \leq
k} \vert v \vert _{n,D} + \lceil v \rceil _{k,\varsigma,D}.
\]}}
\begin{subequations}\label{prob-interior}
\begin{eqnarray}
 Lu \equiv -\ve \triangle u+ au_x +bu=f, \quad (x,y) \in \Omega ,\\
u=0, \  (x,y) \in \partial \Omega ,\quad a > \alpha >0,\  b \geq 0, \quad a,b,f \in C^{5, \gamma }(\bar \Omega  ) .
\end{eqnarray}
\end{subequations}
We   define the {\it inflow boundary} $\partial \Omega _I$ and the {\it outflow boundary} $\partial \Omega _O$ by
\[
\partial \Omega _I:= \{(\phi ,\psi ) \vert n_1 > 0  \} \quad \hbox{and} \quad  \partial \Omega _O:= \{(\phi ,\psi ) \vert   n_1  < 0    \},
\]
where $\vec n=(n_1,n_2)$ is the inward-pointing unit normal. If $\psi '(t_C) =0$ (i.e., $\vec n=(0,\pm 1)$), then this will correspond to a characteristic point $(\phi (t_C) ,\psi (t_C))$ on the boundary.  To exclude the presence of parabolic boundary layers \cite[Chapter 6]{fhmos} and  \cite[Chap. 4, \S 1]{ilin}, we assume that there is only a finite number of isolated characteristic points on the boundary. 
We confine the discussion to  characteristic points where the component $n_1$ changes sign. 

{\bf Assumption 1.}  {\it Assume that there is a finite number of characteristic points on the boundary $\partial \Omega$.
Moreover, at each characteristic point $t=t_C$  assume that there exists  a $\delta >0$ and  a neighbourhood $t_C\in I_C: = (t_0, t _1)$ such that
$\vert \psi '(t_0) \vert = 2\delta = \vert \psi '(t_1) \vert$  and  
\begin{equation}\label{assume0}
 \psi '' (t) >0 \quad \hbox{or} \quad \psi '' (t) <0, \ \forall  t \in I_C. 
\end{equation}
}

We identify three subintervals of $I_C=I^\delta _O\cup I^\delta _C \cup  I^\delta _I.$, associated with each characteristic point $(\phi (t_C) ,\psi (t_C))$:
\begin{subequations}\label{delta-intervals}
\begin{eqnarray}
I^\delta _O(t_C)&:=& \{ t  \in I_C \vert 2\delta \geq \psi ' (t) \geq \delta \} ,\\ I^\delta _I(t_C)&:=& \{ t  \in I_C \vert  -2\delta \leq \psi ' (t) \leq -\delta  \}, \\
I^\delta _C (t_C)&:=& \{ t  \in I_C  \vert  \vert  \psi ' (t) \vert <\delta   \}.
\end{eqnarray} 
\end{subequations}
As the domain is closed, there will be  at least two distinct characteristic points on the boundary $\partial \Omega$.
If the domain has an internal tangent to the boundary at $t=t_C$ and $n_1$ changes sign at this point, 
then  we shall call $(\phi (t_C) ,\psi (t_C))$  an {\it internal characteristic point}.
Otherwise,  if $n_1$ changes sign at $t=t_C$, we call $(\phi (t_C) ,\psi (t_C))$ an {\it external characteristic point} \cite[Chap. 4]{ilin}.

In the local coordinate system, the differential equation transforms into:
\begin{eqnarray*}
-\ve \Bigl( \eta ^{-1} \frac{\partial }{\partial r} \bigl( \eta \frac{\partial \tilde u }{\partial r} \bigr) + 
\zeta  \frac{\partial }{\partial t} \bigl( \zeta \frac{\partial \tilde u }{\partial t} \bigr) \Bigr)+ \tilde a n_1\tilde u_r+ \tilde a \zeta n_2   \tilde u_t +\tilde b u=\tilde f. 
\end{eqnarray*}

 Let us partition the domain $\Omega $ into a finite number of non-overlapping subdomains $\{ D_i \} _{i=1}^n$ such that
\[
D_i \cap D_j = \emptyset , \  i \neq j, \quad \bar \Omega = \cup _{i=1}^n \bar D_i \quad \hbox{and} \quad \Gamma _{i,j} := D_i \cap D_j.
\]
Let $\frac{\partial u }{\partial n^+ _i}$ (and  $\frac{\partial u }{\partial n^-_i}$) denote the outward (and inward) normal derivative of each subdomain $D_i$ and we define the jump in the normal derivative across an interface $\Gamma _{i,j}$ to be
\[
 \ \bigl[ \frac{\partial u }{\partial n} \bigr] _{\Gamma _{i,j}} := \frac{\partial u }{\partial n^+_i}  \vert _{\Gamma _{i,j}} - \frac{\partial u }{\partial n_i^-} \vert _{\Gamma _{i,j}}.
\]
Using the usual proof by contradiction argument (with a separate argument for the interfaces $\Gamma _{i,j}$) we can establish the following 

\begin{theorem}\label{comparison} If $ w,v \in C^0(\bar \Omega ) \cap (\cup _{i=1}^n C^2(D_i ))$  is such that for all $i$, $ L w (x,y)\geq   L v (x,y), \ \forall (x,y) \in  D_i$;  for all $i,j$:  $ \bigl[ \frac{\partial w }{\partial n} \bigr] _{\Gamma _{i,j}} \leq \bigl[ \frac{\partial v }{\partial n} \bigr] _{\Gamma _{i,j}}$ and $ w \geq  v$ on the boundary $\partial \Omega   $, then $ w (x,y) \geq  v(x,y),\ \forall (x,y) \in \bar \Omega$.
\end{theorem}

We next assume  that the  data $a,b, f$   is sufficiently regular  and the boundary $\partial \Omega$ sufficiently smooth so that
$u \in C^{3,\gamma }(\overline{ \Omega }) $. See \cite[pg. 94]{GandT} for definition of smooth domain and boundary. Also see \cite[Theorem 6.14]{GandT} and \cite[Theorem 6.19]{GandT} to justify the following assumption. 

{\bf Assumption 2.}  {\it Assume that   $\Omega $ is a  $C^{3,\gamma }$ domain and the  data $a,b, f \in C^{1,\gamma }(\overline{ \Omega })$,  so that
\begin{equation}\label{assum1}
u \in C^{3,\gamma
}(\overline{ \Omega }) .
\end{equation}
}

As the problem is linear, there is no loss in generality in dealing
with homogeneous boundary data. 
Nevertheless, below we decompose the solution into  regular and layer components, which satisfy a singularly perturbed differential equation with inhomogeneous boundary data.
Hence, we state a result on a priori bounds on the derivatives of the solution of the more general problem:  find $z$ s.t.
\begin{eqnarray}\label{general-problem}
Lz = p(x,y), \quad (x,y) \in  \Omega, \quad z =q(x,y), (x,y) \in  \partial \Omega,
\end{eqnarray}
where the data  and the boundary $\partial \Omega$ are  sufficiently smooth so that $z \in C^{3,\varsigma}(\bar{\Omega})$.

Using stretched variables $\zeta =x/\ve, \iota = y/\ve$ and bounds from \cite{lady}, one can establish the following result (see the argument in \cite[Appendix A]{mos} or \cite[Theorem 3.2]{LS01}).
\begin{lemma}\label{bndsonderivs} Assume $a,b, p \in  C^{1,\varsigma}(\bar{\Omega})$ and $ q\in  C^{3,\varsigma}(\partial {\Omega})$.
The solution $z$ of  \eqref{general-problem}  satisfies 
 \begin{eqnarray*}
||z|| &\leq& C||p||+ ||q||_{\partial \Omega} ,
\\
|{z}|_{1} + \ve  ^\varsigma \lceil {z}\rceil_{1,\varsigma} &\leq&
 C\ve ^{-1} ||z||  + C\Bigl(|{p}|_{0}+   +  \ve ^{\varsigma }\lceil {p}\rceil _{0,\varsigma}
+  \sum _{i=0}^2 \ve ^{i-1}| {q}|_{i}+ \ve ^{1+\varsigma} \lceil  q \rceil _{2,\varsigma , \partial \Omega } ) \Bigr) ,
\end{eqnarray*}
 and for $l = 0,1$
 \begin{eqnarray*}
& | z |_{2+l} + \ve ^\varsigma  \lceil z  \rceil _{l+2,\varsigma} \leq C\ve ^{-(2+l)}||z||   +
 C\Bigl( \sum _{i=0}^l \ve ^{i-(1+l)}|{p}|_{i}+    \sum _{i=0}^l  \ve ^{\varsigma -1}\lceil {p}\rceil _{i,\varsigma}\Bigr) \\
& \qquad +
C\Bigl( \sum_{i= 0}^{l} \ve ^{i-l-2}
|q|_{i} +  \ve ^{\varsigma} \sum_{i= 0}^{l}   \lceil q \rceil _{2+i,\varsigma }\Bigr).
\end{eqnarray*}
\end{lemma}

From these bounds, we have the crude bounds on the solution $u$ of \eqref{prob-interior}
\begin{equation}\label{crude}
\vert u \vert _n  \leq C\ve ^{-n},\quad n =0,1,2,3.
\end{equation}

The solution $u$ of problem (\ref{prob-interior}) can be  decomposed    into the sum
$u=v+w$,
where $v$ is a regular component and $w$ is a boundary
layer function associated with the outflow boundary.
The {\it reduced solution} $v_0$ is defined as the solution of the first order problem
\begin{equation}\label{reduced}
a(v_0)_x+bv_0=f, (x,y) \in \bar \Omega \setminus \partial  \Omega _I, v_0(x,y) =0, \ (x,y) \in \partial   \Omega _I;
\end{equation}
and the first correction $v_1$ is defined as the solution of
\begin{equation}\label{correction}
a(v_1)_x+bv_1=\triangle v_0, (x,y) \in \bar \Omega \setminus \partial  \Omega _I, v_1(x,y) =0, \ (x,y) \in \partial  \Omega _I.
\end{equation}

Our next assumption guarantees that only regular  boundary layers appear near the outflow boundary.

{\bf Assumption 3.}  {\it At each characteristic point $(x_C,y_C)=(\phi (t_C), \psi(t_C)) $,  define the rectangle
\[
Q_\delta (x_C,y_C)
 := \bigl ((-\infty, x_C+\delta] \times [y_C-\delta, y_C +\delta ]\bigr)  \cap \bar \Omega.   
\]
Assume that $a,b, f \in C^{5,\gamma }(\overline{ \Omega })$ and that there exists some $\delta >0$ such that 
\begin{equation}\label{assum2}
f(x,y) \equiv 0 \quad \hbox{for all} \quad (x,y) \in Q_\delta (x_C,y_C), \quad \forall (x_C,y_C).
\end{equation}
}

\begin{remark}
At any external characteristic point $\vec p$, assumption (\ref{assum2}) constrains the data in a $O(\delta ^2)$ neighbourhood of $\vec p$.
If the characteristic point $\vec p$ is an internal characteristic point, then assumption (\ref{assum2}) constrains the data in a $O(\delta )$ neighbourhood of $\vec p$.
\end{remark}

As $v_0,v_1$ are solutions of a first order differential equation (either (\ref{reduced}) or (\ref{correction})), it follows from Assumption 3 that
\[
v_0(x,y) \equiv v_1(x,y)  \equiv 0, \quad \hbox{for all} \quad \vert y - y_C \vert < \delta .
\]

The {\it regular component} $v$  is defined as the solution of the problem: Find $v$ such that
\begin{subequations}
\begin{eqnarray}
Lv=f, \quad (x,y) \in \Omega , \\
v= v_0+\ve v_1,\ (x,y) \in \partial \Omega _O, \ v=u=0, \ (x,y) \in \partial \Omega \setminus \partial \Omega _O.
\end{eqnarray}
\end{subequations}
Note that, if $z=v- (v_0+\ve v_1)$, then $z$ satisfies 
\begin{eqnarray*}
Lz=\ve ^2 \triangle v_1 \ (x,y) \in \Omega , \quad 
z= 0,\ (x,y) \in \partial \Omega .
\end{eqnarray*}
As in \cite[Appendix A]{mos} and since Assumption 3 eliminates any regularity difficulties at the characteristic points, we then have $v_0 \in  C^{5, \gamma }(\bar \Omega  ),\ v_1 \in C^{3, \gamma }(\bar \Omega  ) $ and $z \in  C^{3, \gamma }(\bar \Omega  ) $. Using Lemma \ref{bndsonderivs}, we deduce that
\[
\Vert z \Vert \leq C\ve ^2, \quad \vert z \vert _i  \leq C (1+\ve ^{2-i}), \ 1 \leq i \leq 3.
\]
Hence, we have the following bounds on the regular component
\begin{equation}\label{regular-bnd}
\Vert v \Vert \leq C, \quad \vert v \vert _i + \ve ^{\gamma}\lceil
v \rceil_{i,\gamma} \leq C (1+\ve ^{2-i}), \ 1 \leq i \leq 3.
\end{equation}
Assumption (\ref{assum2})  also  prevents any internal characteristic layers emerging from any  internal characteristic points. 

The {\it boundary layer component} $w$ is the solution of   the problem: Find $w$ such that
\begin{subequations}\label{bdy-layer}
\begin{eqnarray}
Lw=0, \quad (x,y) \in \Omega , \\
w = (u-v) (x,y), \ (x,y) \in \partial \Omega _O \quad w(x,y)=0, \ (x,y) \in \partial \Omega \setminus \partial \Omega _O.
\end{eqnarray}
\end{subequations}
From assumption (\ref{assum1}), it follows that for all characteristic points 
\[
w(x,y) =0, \quad \hbox{if} \ (x,y) \in \partial \Omega _0 \quad \hbox{and} \quad \vert y - y_C \vert < \delta .
\]
In other words, the boundary layer function is not only zero on the inflow boundary, but also on those parts of the outflow boundary 
near the characteristic points. 

 For some fixed $R=O(1)$,
we define the {\it strip}
\begin{eqnarray}\label{strip}
\tilde \Omega _S:= \{ (r, t) \in \tilde \Omega \vert \quad  \psi '(t) >0,\ r \in (0,R) \}, \quad R < \min \{ \delta , \Vert \kappa \Vert _{\partial \Omega _O}  ^{-1}\} ,
\end{eqnarray}
 which is aligned to the outer boundary $\partial \Omega _O$. 
The width $R$ of this strip $\Omega _S$ is further limited in (\ref{maxR}). 
Note that if there are internal characteristic points, then $\Omega _S$ will not be a connected set. The outer boundary $\partial \Omega \cap \partial \Omega _S$ of the strip will have characteristic points as end-points.  At each of  these characteristic points, the  strip $\Omega _S$ will have a vertical  boundary of the form $\{ (x_C,y) \vert \vert y-y_C \vert \leq R\}$. 
If $R < \delta$ then the boundary layer function  $w\equiv 0$ along these vertical  boundaries of  $\Omega _S$, by (\ref{assum1}).

\begin{lemma}\label{main-Lemma} {\it  Assume (\ref{assume0}), (\ref{assum1}), (\ref{assum2}). If $w$ is  the solution of (\ref{bdy-layer}) then within the strip $\Omega _S$ (\ref{strip}),
\begin{eqnarray}\label{bdy-layer-bound}
\vert   \tilde w (r,t) \vert \leq Ce^{-\frac{\alpha \theta r}{ \ve } }+ Ce^{-\frac{\mu \alpha \theta R}{ \ve } },\   0 \leq r \leq R; \ 
\theta  := \min _{\psi ' \geq \delta >0} \vert n_1 \vert, \ \mu < 1; 
\end{eqnarray}
and exterior to  the strip $\vert  w (x,y)\vert \leq Ce^{-\frac{\mu \alpha \theta R}{ \ve } },  \quad (x,y) \in \Omega \setminus \Omega _S.$
}
\end{lemma}
\begin{proof} 
Associated with the neighbourhood $I_C$ (\ref{delta-intervals}) of each  characteristic point, we  construct a cut-off function $\Psi _C(t;2\delta ) \in C^2[0,T]$ such that 
\begin{subequations}
\begin{eqnarray}
0 <\Psi _C(t;2\delta ) < 1,\ t \in I_O^\delta ;\label{cut-off} \\ \Psi _C(t;2\delta )\equiv 1, t \in  \partial \Omega _O \setminus I_C; \ \Psi _C(t;2\delta )\equiv 0, \ t \in \partial \Omega \setminus \partial \Omega _O.
\end{eqnarray}
\end{subequations}
Consider the barrier function
\[
\tilde B(r,t) :=\Psi _C(t;2\delta )\Bigl(  \frac{e^{-\frac{\alpha \theta r}{ \ve } } -e^{-\frac{\alpha \theta R}{ \ve }}}{1 -e^{-\frac{\alpha \theta R}{ \ve }}}\Bigr) , \quad (r,t) \in \overline{ \tilde \Omega _S}.
\]
Observe that $\psi ' (t) \geq \delta,  \ \vert \Psi '_C(t) \vert + \vert \Psi ''_C(t) \vert \leq C, \ t \in I_0^\delta$ and $\Psi _C (t) \equiv 0, \ t \in I_C \setminus I_0^\delta$. 
Then, for $\ve$ sufficiently small and $(r,t) \in \Omega _S$, using the  definition of $\theta$ we have
\begin{eqnarray*}
\tilde L \tilde B =\frac{\alpha \theta }{\ve } \bigl (  \tilde a \frac{\psi ' (t)}{\tau (t)} -\alpha \theta + O(\ve) \bigr ) \tilde B(r,t) > 
\frac{\alpha ^2\theta}{\ve }(-n_1 -\theta  ) \tilde B(r,t) \geq 0. 
\end{eqnarray*}
Observe that $\tilde B = \Psi _C(t;2\delta )$ on the outer boundary $\partial \Omega _O$ and $\tilde B =0 $ on the other three boundaries of the strip. 
This function $\tilde B(r,t)$ is currently only defined on the strip $\bar \Omega  _S$. We extend this function to $\bar \Omega$ 
as follows:
\[
B(x,y) := \Bigl \{ \begin{array}{ll} \tilde B(r,t),\quad (x,y) \in \bar \Omega  _S\\
\quad  0,\qquad  (x,y) \in  \Omega \setminus \bar \Omega  _S
\end{array},
\]
such that  at the inner boundary    $\partial \Omega _S^-:= \partial \Omega _S \cap \Omega$, 
\[
\frac{ \partial  B}{\partial n _s^+} = -\frac{ \partial  \tilde B}{\partial r}= \frac{\Psi _C(t;2\delta )}{\ve} \frac{\alpha \theta e^{-\frac{\alpha \theta R}{ \ve } } }{1 -e^{-\frac{\alpha \theta R}{ \ve }}} > 0 \quad \hbox{and}  \quad \frac{ \partial  B}{\partial n _s^-}=0 .
\]
We   define a second barrier function of the form
\[
B_1(x,y) := \Bigl \{  \begin{array}{ll} \phi (t) - m_x +C_*(R-r), \qquad \quad  (x,y) \in \bar \Omega _S
\\
\  x - m_x,\qquad  (x,y) \in  \Omega \setminus \bar \Omega  _S 
\end{array} ,
\]
where $m_x := \min _{0\leq t \leq T} \phi (t)$. On the strip, for $\ve$ sufficiently small, 
\[
\tilde L \tilde B_1 = \frac{\tilde a}{\tau} (C_*(\psi ')^2 + \zeta (\phi ')^2) + O(\ve) \geq 0.
\]
At the inner boundary  of the strip $\partial \Omega _S^-$, we have
\[
\frac{ \partial  B_1}{\partial n _s^-} = \frac{\psi '}{\tau } \frac{ \partial  B_1 }{\partial x} |_{\partial \Omega ^-_S } =   \frac{\delta }{\Vert \tau \Vert } \quad \hbox{and} \quad \frac{ \partial  B_1}{\partial n _s^+} = C_*:=\frac{\delta }{2\Vert \tau \Vert }.
\] 
We complete the proof by forming the barrier function
\[
B_2(x,y) := B(x,y) +  \frac{2\Vert \tau \Vert}{\delta} \frac{\alpha \theta R}{ \ve }e^{-\frac{\alpha \theta R}{ \ve } }B_1(x,y) .
\]
By the design of this barrier function, we see that
\[
\frac{ \partial B_2 }{\partial n _s^-} \geq \frac{ \partial B_2 }{\partial n _s^+}.
\]
Complete the proof using the inequality 
\[
te^{-t} \leq (1- \mu )^{-1} e ^{-\mu t}, \quad \mu <1, t \geq 0
\]
and the  comparison principle in Theorem \ref{comparison}. 
\end{proof}

 The crude bounds
\begin{equation}\label{layer-crude} 
\vert w \vert _{\bar \Omega , n}  \leq C\ve ^{-n},\quad n =0,1,2,3.
\end{equation}
on the derivatives of the layer function follow from (\ref{crude}).
We  can further decompose the boundary layer function within the strip:
\begin{eqnarray*}
\tilde w = \tilde w_0 + \tilde w_1, \quad (r,t) \in \overline{ \tilde \Omega _S}, \qquad \hbox{where} \\
\tilde L w_0= 
\tilde L w_1 =0, \quad (r,t) \in  \tilde \Omega _S, \qquad  w_0(x_C,y) = w_1(x_C,y) =0; \\
\tilde w_0(0,t) = \tilde w(0,t), \ \tilde w_0(R,t) = 0, \qquad \tilde w_1(0,t) = 0,  \tilde w_1(R,t) = \tilde w(R,t). 
\end{eqnarray*}
By Lemma \ref{main-Lemma} and the maximum principle, $\Vert \tilde w_1 \Vert _{\tilde \Omega _S} \leq  Ce^{-\frac{\alpha \theta R}{ \ve } }$.
Moreover, from the bounds in (\ref{layer-crude}),
\begin{equation}\label{orth-normal}
\Bigl \Vert \frac{\partial ^i
\tilde w_0}{ \partial  r^i} \Bigr\Vert _{ \Omega _S}\leq C(1+\ve ^{-i}), \quad i \leq 3. 
\end{equation}
Consider the function 
$
\tilde z(r,t) = \tilde w_0(r,t) - \Phi(r,t), $
where, for each $t$, $\Phi(r,t)$ is the solution of the problem
\begin{eqnarray*}
-\ve \Phi _{rr} + \frac{\tilde a(0,t) \psi '(t) }{ \tau (t) } \Phi _{r} =0, \ r \in (0,R), \quad
\Phi (0,t)=\tilde w(0,t) , \Phi (R,t)= 0.
\end{eqnarray*}
Then $\tilde z \equiv 0$ on $\partial \Omega _S$ and we can check that 
\[
\vert \tilde L \tilde z(r,t) \vert \leq C\ve ^{-1}  e^{-\frac{\tilde a(0,t) \psi '(t) r}{ \tau (t) \ve } },\ (r,t) \in \Omega _R.
\]
 Applying the arguments from \cite[Theorem 12.4]{mos}, we deduce that
\begin{equation}\label{orth-bnd}
\Bigl \Vert 
\frac{\partial ^j \tilde w_0}{ \partial  t^j} \Bigr\Vert _{ \Omega _S}\leq C(1+\ve ^{1-j}), \quad  j \leq 3 . 
\end{equation}

Let $\bar \Omega \subset \Omega ^*$ , be an extended smooth closed  domain that encloses $\bar \Omega$ 
and is constructed  so that any characteristic point $(x_C,y_C) \in \partial \Omega $ is extended to $(x_C,y^*_C) \in \partial \Omega ^*$ and $\vert y_C^*\vert \geq \vert y_C\vert$. In this way, all the points on the inflow boundary $\partial \Omega _I$ are extended to the inflow boundary  $\partial \Omega _I^*$ and likewise for the outflow boundaries. With each point $\vec p \in \partial \Omega$ on the boundary, we associate $\vec p_* \in \partial \Omega^*$ as the point lying on the outward normal to $\partial \Omega $ and passing through $\vec p$. If the the two boundaries intersect, then $\vec p_* = \vec p$. We define the {\it width of the extension} to be  
\begin{equation}\label{ext-width}
  \delta _E:= \max _{\vec p \in \partial \Omega} \vert \vec p - \vec p_* \vert. 
\end{equation}
 Let the boundary of this extended domain be parameterized  by
\[
\partial \Omega ^*:= \{ (x, y)\vert x =\phi ^* (t), y =\psi ^*(t),\ 0 \leq t \leq T \}, \quad
 \gamma ^* (t) :=(\phi ^*(t), \psi ^* (t));
\]
and  $a^*,f^*$ are  smooth extensions of $a,f$ from $\Omega $ to $\Omega ^*$.
Define $v^*, w^*$ as the solutions of 
\begin{eqnarray*}  
L^*v^*=f^*, \ (x,y) \in \Omega ^*,\quad v^* =0,\    (x,y) \in \partial  \Omega _I^*, v^* = v_0^*+\ve v_1^*, \  (x,y) \in \partial  \Omega ^*_O,\\
L^*w^*=0, \ (x,y) \in \Omega ^*,\quad w^* =0,\    (x,y) \in \partial  \Omega _I^*, w^* =u^*-v^*, \  (x,y)  \in \partial  \Omega _0^*.
\end{eqnarray*}
As $v_0, v_1, v^*_0, v^*_1$ are solutions of first order problems (\ref{reduced}), (\ref{correction}) 
and $\Vert v_0^* \Vert _{\partial \Omega _I} \leq C \delta _E,  \Vert  v^*_1 \Vert _{\partial \Omega _I} \leq C \delta _E$ , then
\[
\Vert (v_0^*+\ve v^*_1) -  (v_0+\ve v_1) \Vert _{\bar \Omega} \leq C \delta _E.
\]
Then, using a comparison principle for the elliptic operator $L$ 
\[
\Vert v^* - v \Vert _{\bar \Omega} \leq C \delta _E (x+C_1) \leq C  \delta _E.
\]
Hence, on any subdomain $D \subset \bar  \Omega$, we have 
\[
\Vert v^* - u \Vert _{D} \leq  \Vert v^* - v - w \Vert _{D} \leq C  \delta _E + \Vert w  \Vert _{D}.
\]
We define an  approximation $u_1$ to $u$ on the strip $\Omega _R$  as  the solution of
\[
\tilde L \tilde u_1=\tilde f, (r,t) \in \tilde \Omega _R, \tilde u_1(0,t)=\tilde u(0,t)0, \tilde u_1(R,t)=\tilde v^*(R,t),\ t \in I_O
\]
where $I_O$ denotes the subinterval of $[0,T]$ where $\psi '(t) >0$. Outside the strip $u_1:=v^*, (x,y) \in \Omega \setminus \Omega _R$. 
Then
\[
\Vert u -u_1 \Vert _\Omega \leq C\delta _E+ \Vert w  \Vert _{\Omega \setminus \Omega _R}.
\]
In the next section, we describe a numerical method which initially generates an approximation $\bar U_0$  to $v^*$ across $\bar \Omega$ and, using  $\bar U_0$ along  the inner boundary $\Omega^-_S$,  it corrects this initial approximation within the strip using a piecewise-uniform Shishkin mesh. 

\section{Numerical method}

Enclose  the domain with a rectangle $\bar Q_E:= [m_x,M_x] \times [m_y,M_y],$ where  
\begin{subequations}\label{Outer-rectangle}
\begin{eqnarray}
m_x \leq  \min _{t\in [0,T]}  \phi (t)  , \quad M_x \geq  \max _{t\in [0,T]} \phi (t)   , \quad  L_x := M_x-m_x;\\
m_y \leq \min _{t\in [0,T]}  \psi (t)  , \quad M_y \geq  \max _{t\in [0,T]}  \psi (t)  ,\quad L_y := M_y-m_y.
\end{eqnarray}
\end{subequations}
Set $u(x,y) \equiv 0, \ (x,y) \in \bar Q_E \setminus \Omega$ for points outside $\Omega$, then solve the problem (\ref{prob-interior})  using upwinding on a uniform rectangular mesh
 \[
\bar Q^{N} := \{ (x_i,y_j) \vert x_i=m_x+i\frac{L_x}{N}, y_j=m_y+ j\frac{L_y}{N} \}_{i,j=0}^N. 
\]
That is, find $U_0$ such that
{\footnote{The
finite difference operators are defined as
\begin{eqnarray*}
D^+_xZ(x_i,y_j) :=\frac{Z (x_{i+1},y _j)-Z(x_i,y _j)}{x_{i+1}-x_i},\quad  D_x^-Z(x_i,y _j) :=\frac{Z(x_{i},y _j)-Z (x_{i-1},y _j)}{x_{i}-x_{i-1}}; \\
2(bD_x^\pm)Z := (b-\vert b \vert) D_x^+Z + (b+\vert b \vert) D_x^-Z;\quad \delta^2_x Z(x_i,y _j) :=\frac{D^+_xZ(x_{i},y_j)-D_r^-Z(x_i,y _j)}{(x_{i+1}-x_{i-1})/2}.\end{eqnarray*}
}}
\begin{subequations}\label{probA}
\begin{eqnarray}
 L^N U_0(x_i,y_j)&=& f(x_i,y_j), \quad (x_i,y_j) \in \Omega \cap Q^{N};\\
 U_0(x_i,y_j)&=&0, \quad  (x_i,y_j) \in \bar Q^{N} \setminus \Omega; \\
\hbox{where} \qquad L^N &:=&  -\ve (\delta _x^2 +\delta ^2_y) + aD_x^-+bI.
\end{eqnarray}
\end{subequations}
Use bilinear interpolation to form an initial global approximation $\bar U_0$ to $u$, defined as
\[
\bar U_0 (x,y) = \sum _{i,j=1}^{N-1} U_0(x_i,y_j) \phi _i(x) \phi ^j(y), \quad (x,y) \in \bar \Omega,
\]
where  $\phi _i(x)$ ($\phi ^j(y)$) is the standard piecewise linear hat functions centered at $x=x_i (y=y_j)$. Since we are using upwinding, the numerical solution will be stable; in the sense that the operator $L^N$ satisfies a discrete comparison principle of the form: If $Z(x_i,t_j) \geq 0, (x_i,t_j) \in \bar Q^{N} \setminus \Omega $ and $ L^NZ(x_i,t_j) \geq 0, (x_i,t_j) \in  Q^{N}$, then $ Z(x_i,t_j) \geq 0, (x_i,t_j) \in \bar Q^{N}$. However,  the layers at the outflow will be smeared and $\bar U_0$ will not be accurate in the boundary layer region. 

 We  correct this approximation in the strip (\ref{strip}), where the width of the strip is such that
\begin{equation}\label{maxR}
R < \min \{ M_x  , \vert m_y \vert , M_y , \frac{1}{\Vert \kappa \Vert _{ \partial \Omega _O}} \}.
\end{equation}
On the strip, we will use  the following upwinded finite difference operator
\begin{eqnarray*}
\tilde L^N\tilde U_{i,j}:= \Bigl( -\ve \eta _{i,j}^{-1} \bar h_i^{-1}\bigl(\eta _{i+1/2,j} D^+_r-\eta _{i-1/2,j} D^-_r  \bigr) + (a n_1)_{i,j} D_r ^\pm  \\
-\ve \zeta _{i,j} \bar k_j^{-1}\bigl( \zeta _{i,j+1/2} D^+_t - \zeta _{i,j-1/2}  D^-_t\bigr)  
 + (an_2\zeta )_{i,j} D_t ^\pm +b _{i,j} \Bigr) \tilde U_{i,j} ;
\end{eqnarray*}
where
\[
\tilde U_{i,j}:=\tilde U_1(r_i,t_j) \quad \hbox{and} \quad \eta _{i+1/2,j} := \eta (\frac{r_i+r_{i+1}}{2}, t_j). 
\]
The operator $\tilde L^N$ again satisfies a discrete comparison principle on the strip. 
Over $\Omega _S$, we generate  a Shishkin mesh \cite{fhmos} in the $r$ coordinate and  a uniform mesh in the $t$ direction. The interval $[0,R]= [0, \sigma] \cup [\sigma , R]$ and  the Shishkin transition point \cite{fhmos} is taken to be
\begin{equation}\label{transition-point}
 \sigma := \min \{ \frac{R}{2} ,  C_*\frac{\ve }{\alpha} \ln N \}, \ C_* >  \frac{1}{\theta},\quad \theta =   \min _{\psi ' \geq \delta >0} \vert n_1 \vert .
\end{equation}

\begin{remark} In the neighbourhood of each characteristic point $t=t_C$, consider the parabola 
\[
y-y_C= a(x-x_C)^2, \quad \vert a \vert = \vert \kappa (t_c) \vert, \quad \hbox{on which}
\]
\[
\vert n_1(t) \vert = \frac{2\vert \kappa (t_c)  t\vert }{\sqrt{1+(2 \kappa (t_c)  t)^2}}.
\]
If we choose the parameter $\delta$ in Assumption 3 such that $  \vert \kappa (t_c) \vert ^{-1} > \delta \geq  (m \vert \kappa (t_c) \vert )^{-1},\ m > 1$ and if the boundary $\partial \Omega$  coincided with this parabola then $\theta \geq \frac{2}{\sqrt{4+m}}. $
Based on this observation, we  shall take $C_*=2$ in (\ref{transition-point}) in our numerical experiments. 
\end{remark}
The mesh on the strip will be denoted by $\tilde \Omega ^N_S$ and the mesh points on the boundary of the strip will be denoted by $\partial \tilde \Omega ^N_S$. 
 Solve for the corrected approximation in the transformed co-ordinates
\begin{subequations}\label{probB}
\begin{eqnarray}
\tilde L ^N\tilde U_1(r_i,t_j) =\tilde f(r_i,t_j) , \quad (r_i,t_j) \in \tilde \Omega ^N_S, \\
\tilde U_1(r_i,t_j) = 0, \ (r_i,t_j)  \in \partial \tilde \Omega ^N_S \cap \partial \Omega _O, \ \tilde U_1= \tilde U_0, \  (r_i,t_j)  \in \partial \Omega ^N_S \setminus  \partial \Omega _O.
\end{eqnarray}
\end{subequations}
Use bilinear interpolation within the strip to form $\bar U_1$. Our corrected  numerical approximation $\bar U$  is  defined by 
\begin{equation}\label{discrete -solution}
\bar U (x,y)= \Bigl\{ \begin{array}{ll} \bar U_1(x,y),\ (x,y) \in \bar \Omega _S, \\
\bar U_0(x,y), (x,y) \in \Omega \setminus \Omega _S\end{array}.
\end{equation}

\begin{theorem}\label{Main-Result}  Assume (\ref{assume0}), (\ref{assum1}), (\ref{assum2}) and that the strip width $R < \delta$.  
If $U$ is the corrected numerical solution (\ref{discrete -solution})  and $u$ is the continuous solution of  (\ref{prob-interior}) then
\[
\Vert \bar U -u \Vert  \leq  CN^{-1} (\ln N)^2.
\]
\end{theorem}

\begin{proof}
In the first phase of the numerical algorithm, we solve on the rectangular mesh $\bar Q^{N}$. For each vertical height $y=y_j$, we identify the external mesh points $S^j := \bar Q^{N}\setminus \bar \Omega$ and the edge co-ordinates
\[
x_{0,j} := \max _{ S^j}  \{ x_i  <0\} \quad \hbox{and} \quad x_{N,j} := \min _{ S^j}  \{ x_i >0\}.
\]
A smooth curve $\partial \Omega ^*$ can be created to pass through these edge  mesh points $\{(x_{0,j}, y_j), (x_{N,j}, y_j)\} _{j=0}^N$, which will define the boundary of an extension $\Omega ^*$ of the domain $\Omega$. Note that the width of this extension (\ref{ext-width}) is such that
$
\delta _E \leq C N^{-1}.
$
Decompose the initial approximation $U_0=V_0+W_0$, where
\begin{eqnarray*}
L^N V_0 =f(x_i,y_j),  \ (x_i,y_j) \in \Omega, \ V_0=v^*, \quad (x_i,y_j) \in \bar Q^{N} \setminus \Omega ,\\
L^N W_0 =0,  \ (x_i,y_j) \in \Omega, \quad  W_0=w^*, \quad (x_i,y_j) \in \bar Q^{N} \setminus \Omega.
\end{eqnarray*}
Let us first bound the error in the regular component. The truncation error and the interpolation error on the boundary yield
\begin{eqnarray*}
\vert L^N (V_0 -v)(x_i,y_j)  \vert = \vert (L^Nv -Lv)  (x_i,y_j) \vert \leq CN^{-1},  \quad  (x_i,y_j) \in \Omega, \\  
\vert (V_0 -v) (x_i,y_j) \vert =\vert (v^*-v) (x_i,y_j) \vert \leq CN^{-1},
 \quad (x_i,y_j) \in \bar Q^{N} \setminus \Omega. 
\end{eqnarray*}
Hence,
\begin{equation}\label{bnd-reg}
\Vert \bar V_0 -v \Vert \leq CN^{-1}.
\end{equation} 
If $\ve$ is sufficently large such that $\ve \ln N >C$, then using the derivative bounds (\ref{crude}) we can deduce that
\begin{equation}\label{bnd-classical}
\Vert \bar U_0 -u \Vert \leq CN^{-1} (\ln N)^2, \quad \hbox{if} \quad \ve \ln N > C.
\end{equation} 
In the other case where $\ve \ln N \leq C$, observe that
\begin{eqnarray*}
\vert W_0(x_i,y_j) \vert \leq C (1+\frac{\alpha H}{\ve})^{i-N} \leq C (1+\frac{\alpha L_x}{N\ve})^{-pN}, \qquad x_i \leq (1-p)N, \ p <1.
\end{eqnarray*}
 As in \cite[Lemma 5.1]{mos}, we  have the inequality{\footnote{Take the natural logarithm of both sides and use  $\ln (1+t) \geq (t(1-t/2), \ \forall t \geq 0$. }},
\[
(1+\frac{q\ln N}{N})^{-pN} \leq  N^{-pq }, \quad  \forall N \geq 1,\ p,q >0.
\]
If $\ve \ln N \leq pq$,  then for all $x_i \leq (1-p)N, p <1$,
\[
\vert W_0(x_i,y_j) \vert \leq C (1+\frac{q}{N\ve})^{-pN} \leq (1+\frac{\ln N}{pN})^{-pN} \leq CN^{-1}.
\]
Then, in the outer domain $\Omega _O:= \Omega \setminus \Omega _S$, for  $\ve \ln N \leq pq$,  we have the bound
\[
\Vert \bar U_0 -u\Vert _{ \Omega_0} \leq \Vert \bar V_0 -v\Vert _{ \Omega_0} + \Vert w \Vert _{ \Omega_0}+ \Vert \bar W_0 \Vert _{ \Omega_0} \leq CN^{-1}.
\]
Hence, in all cases, outside the strip
\[
\Vert \bar U_0 -u\Vert _{ \Omega_0} \leq CN^{-1} (\ln N)^2.
\]

We now examine the error $\Vert \bar U_1 -u\Vert $ on the strip.
Over the interval $I_i:=[x_{i-1},x_{i+1}]$, with  $H_i :=  \max \{h_i, h_{i+1} \} $, we have the truncation error bound
 {\footnote{
\begin{eqnarray*}
&& \frac{h_i}{\bar h_i} D_x^-(a(x_{i+1/2}) D^+_xu(x_i)) = \frac{1}{\bar h_i} \Bigl( \frac{a(x_{i+1/2})}{h_{i+1}} \int _{s=x_i} ^{x_{i+1}} u'(s) \ ds -  \frac{a(x_{i-1/2})}{h_{i}} \int ^{x_i} _{s=x_{i-1}} u'(s) \ ds \Bigr) \\
&=& \frac{1}{\bar h_i} \Bigl( \frac{a(x_{i+1/2})-a(x_i)}{h_{i+1}} \int _{s=x_i} ^{x_{i+1}}  \int _{t=x_i} ^{s} u''(t) \ dt \ ds -  \frac{a(x_{i-1/2})-a(x_i)}{h_{i}} \int ^{x_i} _{s=x_{i-1}} \int _{t=x_i} ^{s} u''(t) \ dt \ ds \Bigr) \\
&+&  \frac{u'(x_i)}{\bar h_i} \int _{s=x_{i-1/2}} ^{x_{i+1/2}}  \int _{t=x_i} ^{s}a''(t) \ dt \ ds +(a'u')(x_i) + a(x_i) \delta ^2_x u(x_i)\\ \\ 
&=& \Bigl(  \frac{a(x_{i+1/2})-a(x_i)}{\bar h_i h_{i+1}} \int _{s=x_i} ^{x_{i+1}}  \int _{t=x_i} ^{s} \int _{w=x_i} ^{t}  -  \frac{a(x_{i-1/2})-a(x_i)}{\bar h_i h_{i}} \int ^{x_i} _{s=x_{i-1}} \int _{t=s} ^{x_i} \int _{w=x_i} ^{t}\Bigr)   u'''(w) \ dw\ dt \ ds \\
&+&  \frac{u''(x_i)}{\bar h_i}\bigl(\frac{h_{i+1}}{2} \int _{s=x_{i}} ^{x_{i+1/2}}  \int _{t=x_{i}} ^{s} a''(t) \ dt \  ds - \frac{h_{i}}{2} \int _{s=x_{i-1/2}} ^{x_{i}}  \int _{t=x_{i}} ^{s} a''(t) \ dt \  ds\bigr) \\
&+&  \frac{u'(x_i)}{\bar h_i} \int _{s=x_{i-1/2}} ^{x_{i+1/2}}  \int _{t=x_i} ^{s} \int _{w=x_i} ^{t} a'''(t) \ dw \ dt \ ds\ 
+  \frac{(u'a''+2u''a')(x_i)}{4}(h_{i+1}-h_i)+(a'u')(x_i) + a(x_i) \delta ^2_x u(x_i)
\end{eqnarray*}
}}
\begin{eqnarray*}
\bigl \vert (au_x)_x(x_i)- \frac{h_i}{\bar h_i} D_x^-(a(x_{i+1/2}) D^+_xu(x_i) \bigr \vert \leq 
C \vert h_{i+1} -h_i\vert \sum _{n=1}^3  \bigl \vert \frac{\partial ^{3-n} a}{\partial x^{3-n}} (x_i)  \bigr\vert  \cdot \bigl\vert \frac{\partial ^n u}{\partial x^n} (x_i)  \bigr \vert  \\
+ CH_i^2 \sum _{n=1}^2 \Vert \frac{\partial ^{4-n} a}{\partial x^{4-n}}  \Vert _{I_i} \cdot \vert \frac{\partial ^n u}{\partial x^n} (x_i)  \vert +CH_i\min \bigl \{ 
\Vert \frac{\partial ^3 u}{\partial x^3}  \Vert _{I_i},  H_i\sum _{n=3}^4 \Vert \frac{\partial ^n u}{\partial x^n}  \Vert _{I_i} \bigr \}.
\end{eqnarray*}
By the choice of transition point (\ref{transition-point}), we can identify the following barrier function
\begin{equation}\label{Barrier}
\tilde E(r_i) := \frac{\Pi _{n=1}^i (1+\frac{\alpha \theta h_i}{\ve})}{\Pi _{n=1}^N (1+\frac{\alpha \theta h_i}{\ve})},\quad \hbox{so that} \quad \tilde L \tilde E (r_i) \geq 0, \ 0<r_i <R.
\end{equation}
The proof is completed using the arguments in \cite[Chapter 13]{mos}, the barrier function (\ref{Barrier}),  Lemma \ref{bdy-layer-bound}, the above truncation error bounds coupled  with the bounds on the derivatives in (\ref{crude}) and (\ref{orth-bnd}).
\end{proof}

\begin{remark}
Recall that throughout this paper, we have assumed that the boundary $\partial \Omega$  is smooth. This assumption is implicitly used in the proof of Theorem \ref{Main-Result} as all the  derivatives $\frac{\partial ^i \zeta }{\partial t^i}, \frac{\partial ^i \eta }{\partial t^i}, \ 0\leq i \leq 3$ are assumed to be bounded within the strip. 
Hence, for the error analysis, the smoothness of the outflow boundary $\partial \Omega _O$  is of particular importance in the error analysis. 
\end{remark}

\section{Numerical results}

The following set of  examples are all related to the following: Let $\beta >0 $ and
\begin{eqnarray*}
\partial \Omega := \{ (\rho \cos t, \rho \sin t) \}, \ \hbox{or} \  \partial \Omega := \{ (\rho \sin t, \rho \cos t) \},  \\
\ \rho (t) := \beta \pm  t^n\sin ^2 t  \quad \hbox{or} \ \rho (t) := \beta \pm  t^n\cos ^2 t\quad  \hbox{with}\quad  n=0,1,2 . 
\end{eqnarray*}
 Once $\phi, \psi$ are smooth functions, the level of smoothness of the boundary $\partial \Omega$ will be determined by identifying the value of $n$ for which 
$\phi ^{(i)} (0)= \phi ^{(i)} (2\pi), \ \psi ^{(i)} (0)= \psi ^{(i)} (2\pi), 0 \leq i \leq n$. (See Example 2 below). 

To numerically estimate the order of convergence of the numerical method as it is applied to several test problems, we use the double-mesh method \cite{fhmos}. 
Denote the mesh points on the rectangular grid $\bar Q^N$ , which lie within the domain $\bar \Omega$ by $\bar Q_I^N:= \bar Q^N \cap \bar \Omega$. 
For each  particular value of $\ve \in R_\ve :=\{ 2^{-i}, i=0,1,2,\ldots 20 \}$ and $N \in R_N: =\{ 2^{-j}, j=3,4,5\ldots 10 \}$, let  $U^N$ be the computed solutions numerical solution (\ref{discrete -solution}), where $N$ denotes the number of mesh elements used in each co-ordinate direction within the rectangle $\bar Q^N$ and within the strip $\bar \Omega _S$.  Define the maximum local two-mesh global differences $ D^N_\ve$ and the  parameter-uniform two-mesh global differences $D^N$ by
\[
D^N_\ve:= \max \{  \Vert \bar U_0^N-\bar U_0^{2N}\Vert _{(Q_I^N \cup Q_I^{2N}) \setminus \Omega _S}, \Vert \bar U_1^N-\bar U_1^{2N}\Vert _{\tilde \Omega ^N_S \cup \tilde \Omega_S^{2N}} \} ;   D^N:= \max _{ \ve \in R_\ve} D^N_\ve .
\]
 Then, for any particular value of $\ve $ and $N$, the local orders of global convergence are denoted by $\bar p ^{N}_\ve$ and, for any particular value  of $N$ and {\it all values of $\ve $}, the  {\bf parameter-uniform} global orders of   convergence $\bar p ^N$ are defined, respectively,  by 
\[
 \bar p^N_\ve:=  \log_2\left (\frac{D^N_\ve}{D^{2N}_\ve} \right) \quad \hbox{and} \quad  \bar p^N :=  \log_2\left (\frac{D^N}{D^{2N}} \right).
\]
In implementing the numerical method, we will apply the method to problems which do not satisy the assumptions (\ref{assume0}), (\ref{assum1}), (\ref{assum2}). Hence, unless otherwise indicated, we simply take $R=0.1$ (for the width of the strip) and $C_*=2 $ (in (\ref{transition-point})). \vskip0.5cm

\begin{figure}[h!]
\centering
\resizebox{\linewidth}{!}{
	\begin{subfigure}[Domain $\Omega _1$ with $\beta =0.5$]{
		\includegraphics[width=7cm, scale =5]{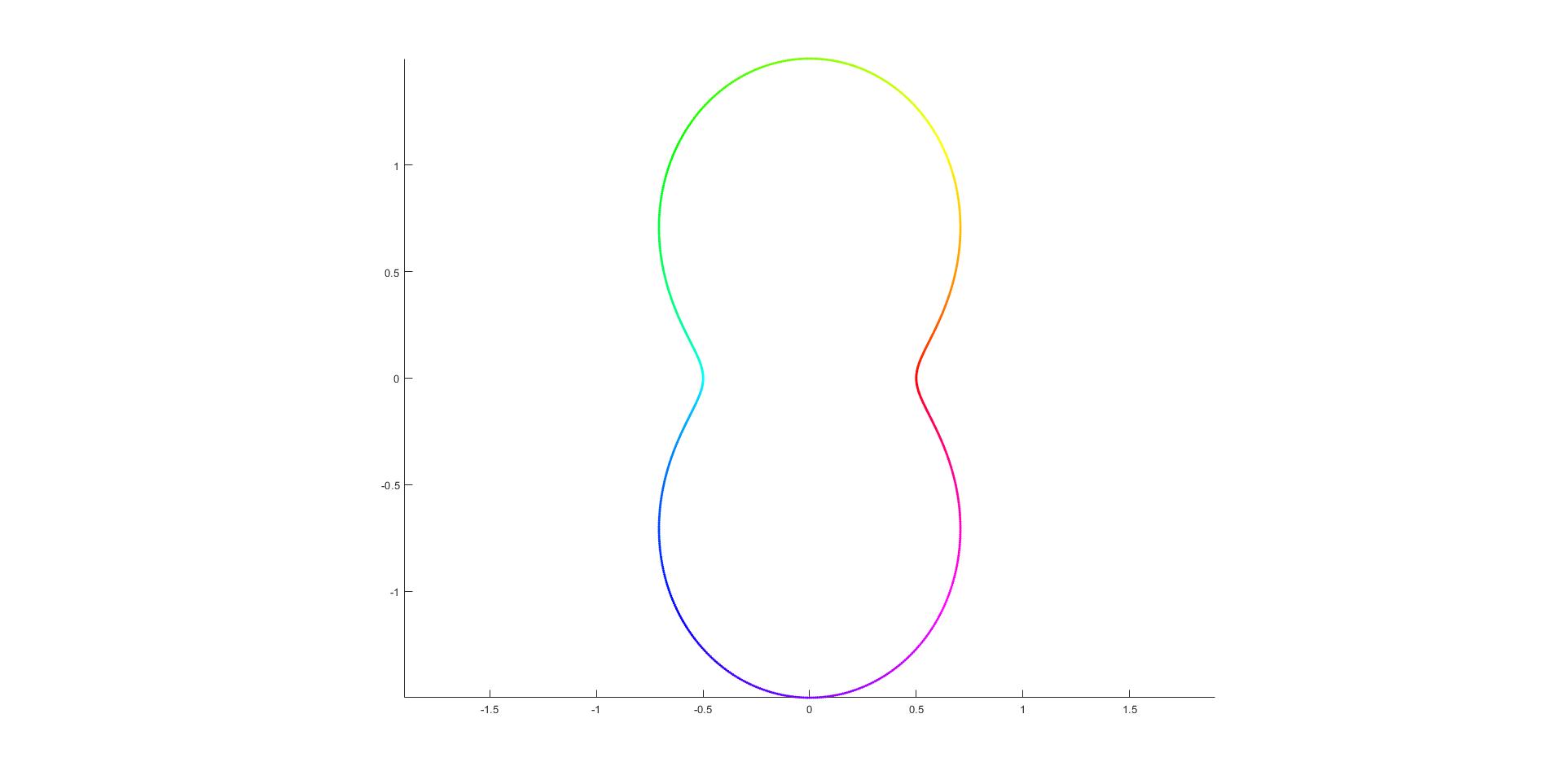} 
		}
    \end{subfigure}
\begin{subfigure}[Computed solution for $\ve =10^{-3}$]{
		\includegraphics[width=7cm, scale =5]{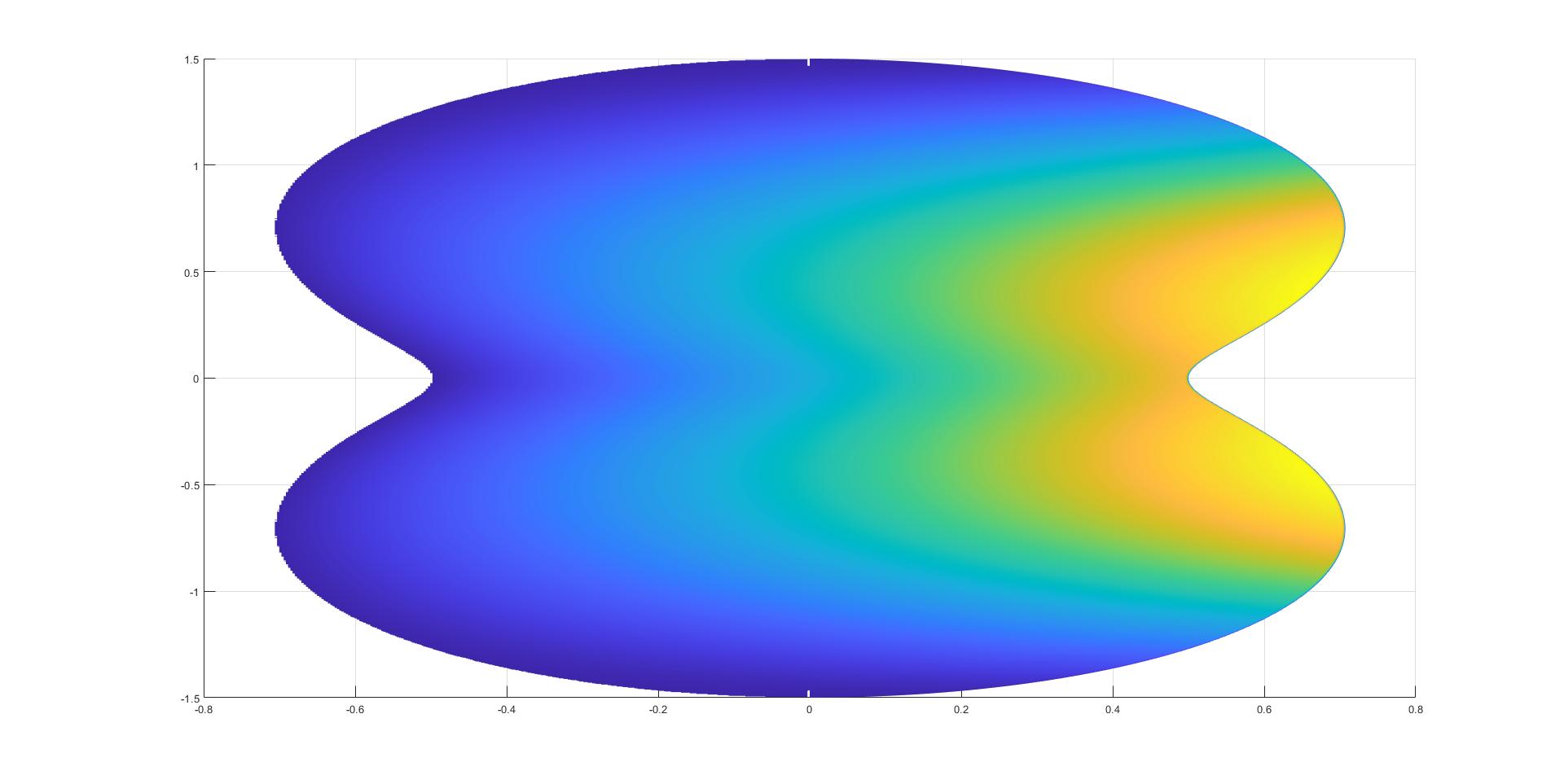} 
		}
	\end{subfigure}
}
	\caption{Computed solution of problem (\ref{Test2}) for $\ve =10^{-3}$ on the domain $\Omega _1$ with $\beta =0.5$}
	\label{fig:ex1-solution}
 \end{figure}

\noindent {\bf Example 1} Consider $\partial \Omega _1:=\{(\phi(t), \psi (t)) \} = \{ (\rho \cos t, \rho \sin t) \}$,
with $\rho := \beta + \sin^2 t$. The domain is displayed in Figure \ref{fig:ex1-solution} for the case of $\beta =0.5$. The outflow is all points on the boundary where $x >0$ and the inflow is for the boundary points where $x<0$.  There are two external characteristic points at $(0, \pm (1+ \beta ))$ and no internal characteristic points. The curvature at the two characteristic points is $\kappa = (3+\beta)(1+\beta )^{-2}$ and the upper bound in (\ref{maxR}) is $R < \frac{1}{6}$.

In Table \ref{Tab3} we present orders   for  the following test problem:
 \begin{equation}\label{Test2}
 -\ve \triangle u+ u_x +u=(1+\beta)^2-y^2, \  (x,y) \in \Omega _1, \quad u=0, \ (x,y) \in \partial \Omega _1
\end{equation}
  A  sample computed solution is displayed in Figure \ref{fig:ex1-solution}.

 \begin{table}[ht]
\centering\small
\begin{tabular}{|r|rrrrrrr|}\hline 
$\varepsilon| N$&$N=8$& 16&32&64&128&256&512\\ \hline 
1.000000 &   0.5137 &   0.5011 &   1.4917  &   0.2148  &   1.7443 &   0.0500 &   1.3869\\ 
0.500000 &   0.1955 &   0.6965 &   1.3812  &   0.2949  &   1.7130 &   0.0752 &   1.3852\\ 
0.250000 &  -0.1489 &   0.8925 &   1.0994  &   0.5916  &   1.2696 &   0.5234 &   1.3834\\ 
0.125000 &  -0.1181 &   0.9350 &   1.0108  &   0.7912  &   1.0026 &   0.8120 &   1.3808\\ 
0.062500 &   0.1995 &   0.8365 &   0.3256  &   1.4904  &   0.7879 &   1.0620 &   1.2017\\ 
0.031250 &   0.4192 &   0.5178 &  -0.0990  &   1.5718  &   0.6644 &   1.0636 &   1.0025\\ 
0.015625 &   0.5074 &   0.5854 &  -0.1666  &   1.5261  &   0.9509 &   1.1332 &   1.0657\\ 
0.007813 &   0.5821 &   0.8212 &   0.0382  &   1.8144  &   1.2916 &   0.6715 &   1.3201\\ 
0.003906 &   0.6319 &   0.8092 &   0.7434  &   1.1094  &   1.2111 &   0.7548 &   1.2573\\ 
0.001953 &   0.6607 &   0.8019 &   0.9146  &   0.9390  &   1.0947 &   0.8686 &   1.1738\\ 
0.000977 &   0.6763 &   0.7984 &   0.9101  &   0.9444  &   1.0248 &   0.9377 &   1.0889\\ 
0.000488 &   0.6843 &   0.7968 &   0.9093  &   0.9455  &   0.9887 &   0.9737 &   1.0304\\ 
0.000244 &   0.6884 &   0.7961 &   0.9093  &   0.9455  &   0.9796 &   0.9558 &   1.0201\\ 
0.000122 &   0.6905 &   0.7958 &   0.9095  &   0.9453  &   0.9830 &   0.9489 &   0.9877\\ 
0.000061 &   0.6915 &   0.7956 &   0.9096  &   0.9452  &   0.9823 &   0.9643 &   0.9623\\ 
0.000031 &   0.6920 &   0.7955 &   0.9097  &   0.9452  &   0.9806 &   0.9818 &   0.9536\\ 
0.000015 &   0.6923 &   0.7955 &   0.9097  &   0.9451  &   0.9787 &   0.9837 &   0.9731\\ 
0.000008 &   0.6924 &   0.7955 &   0.9098  &   0.9451  &   0.9778 &   0.9846 &   0.9941\\ 
0.000004 &   0.6925 &   0.7955 &   0.9098  &   0.9451  &   0.9773 &   0.9851 &   0.9940\\ 
0.000002 &   0.6925 &   0.7954 &   0.9098  &   0.9451  &   0.9769 &   0.9855 &   0.9604\\ 
0.000001 &   0.6926 &   0.7954 &   0.9098  &   0.9451  &   0.9766 &   0.9780 &   0.9448\\ \hline 
$\bar p ^N$ &   0.6558 &   0.6238 &  -0.1517  &   1.5261 &   0.7629 &   1.0636  &   1.0025\\  \hline
\end{tabular}
\caption{Computed double-mesh global orders of convergence $\bar p_\ve ^N$ for the corrected approximations $\bar U$ 
for problem (\ref{Test2}) and  $\beta =0.5$.}
\label{Tab3}

\end{table}

\begin{remark} In determing the double-mesh global orders of convergence, we note that we can underestimate the orders due to a potential overestimate of the maximum two mesh differences. This overestimate is caused by the fact that the interpolant $\bar U_0$
over the rectangular grid, will use mesh points lying within the strip $\bar \Omega _S$. This overspill is visible in Figure \ref{overspill}. An overestimate can occur when the maximum two mesh differences are located near the inner boundary of the strip.
 \begin{figure}[h!]
\centering
\includegraphics[width=5cm]{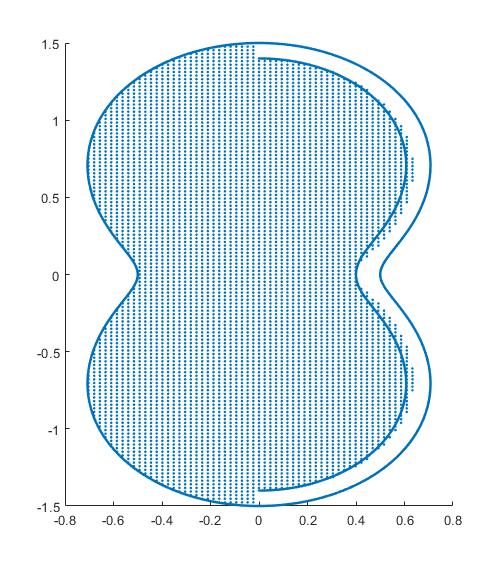}
	\caption{Example 1: Location of mesh points used in the calculation of $\bar U_0(x,y),\ (x,y) \in \Omega \setminus \Omega _S$,  with $N=128$ on the rectangular grid}
	\label{overspill}
 \end{figure}

\end{remark}

\begin{figure}[h!]
\centering
\resizebox{\linewidth}{!}{
	\begin{subfigure}[Domain $\Omega _2$ with $\beta =0.5$]{
		\includegraphics[width=7cm, scale =5]{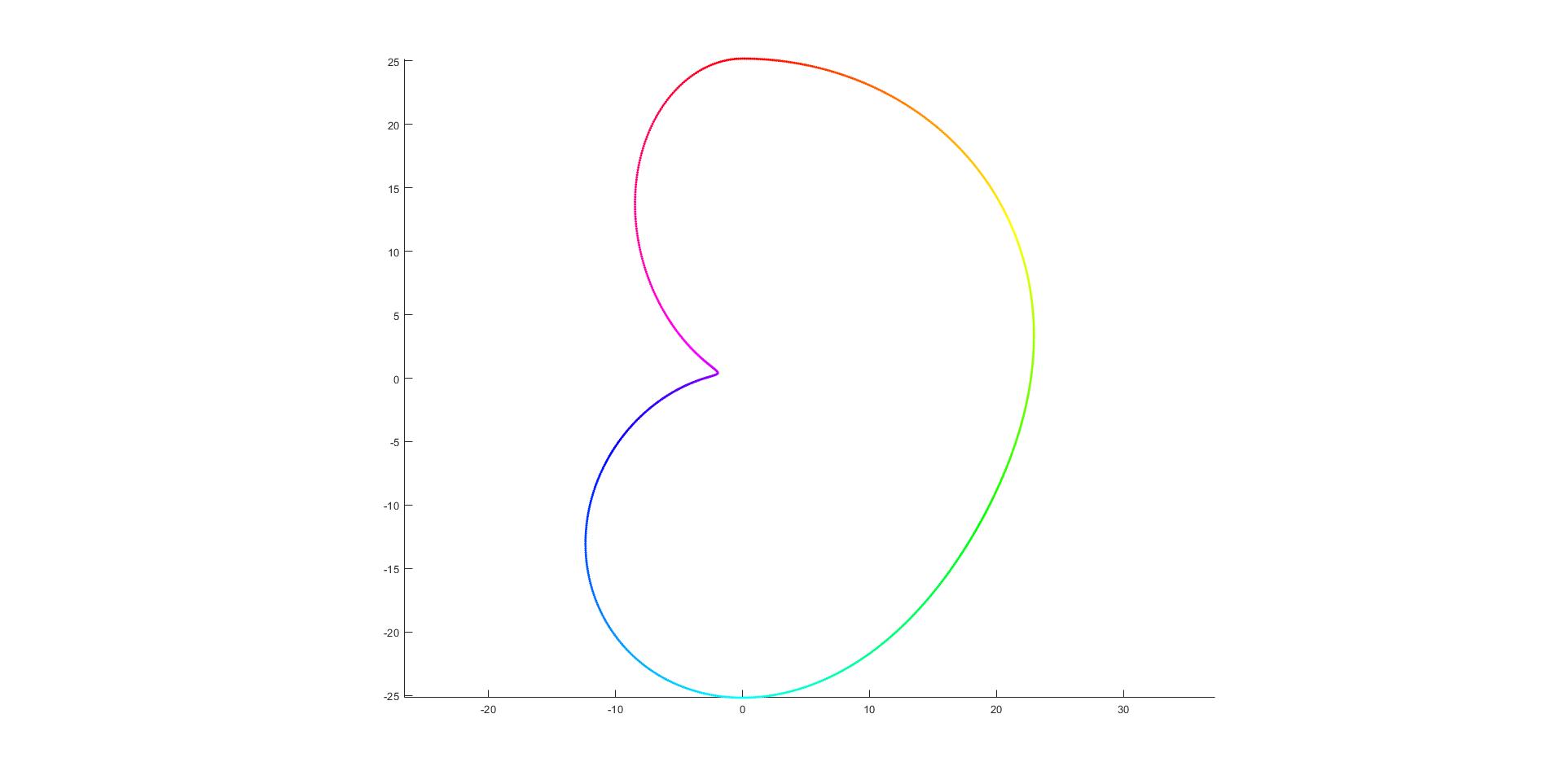} 
		}
    \end{subfigure}
\begin{subfigure}[Computed solution for $\ve =2^{-10}$]{
		\includegraphics[width=7cm, scale =5]{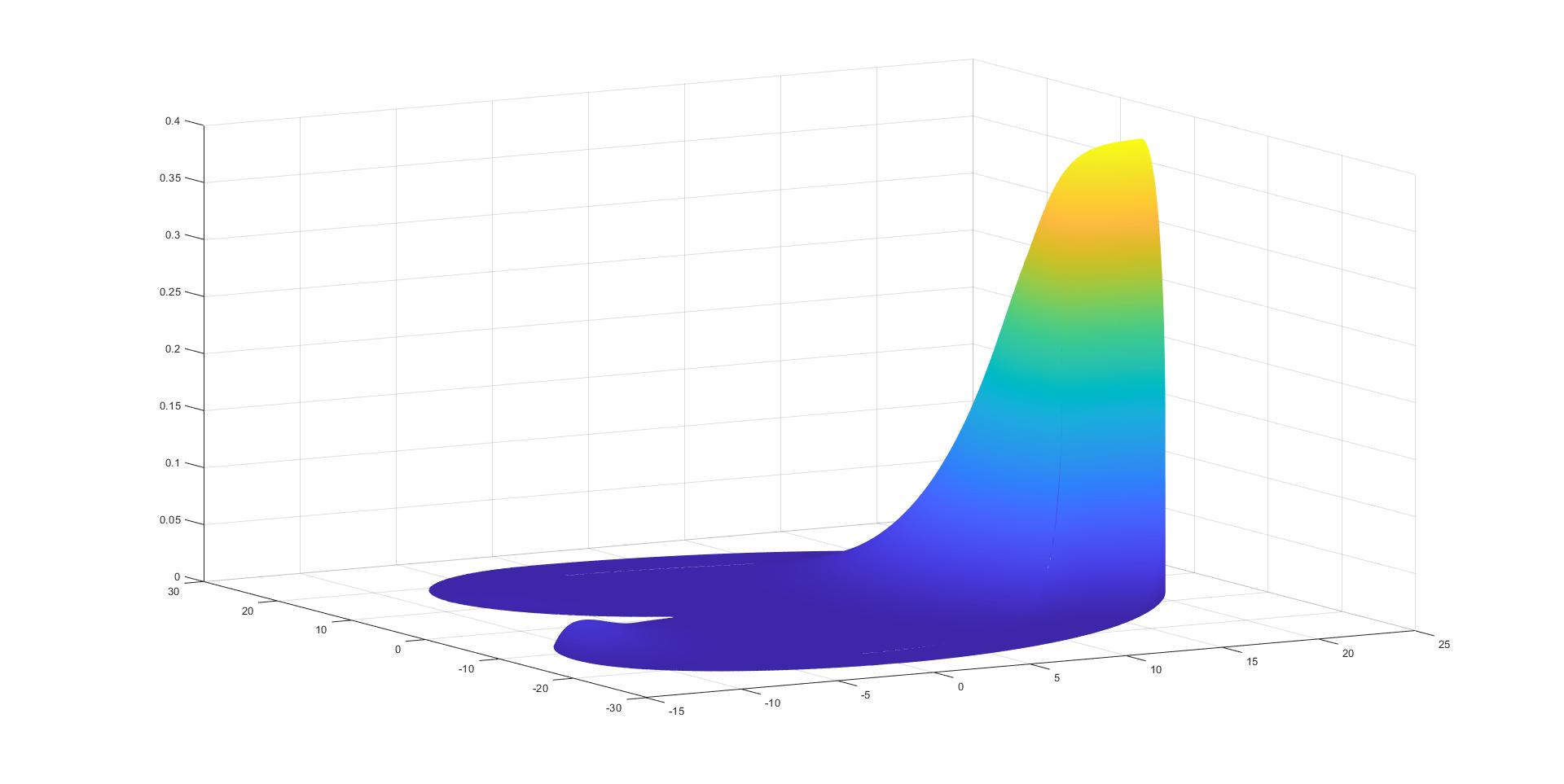} 
		}
	\end{subfigure}
}
	\caption{Computed solution of problem (\ref{Test3}) for $\ve =2^{-10}$ on the domain $\Omega _2$ with $\beta =0.5$}
	\label{fig:ex2-solution}
 \end{figure}

\noindent {\bf Example 2} The boundary is $\partial \Omega _2:=\{ (\rho \sin t, \rho \cos t) \}$, with $\rho := 2.5*\pi ^2 +\beta  -t^2 \sin^2 t$.
In this case the orientation of the curve is clockwise as the parameter $t$ increases and the domain is not symmetrical.  This domain does not have a smooth boundary at $t=0$, where there is a jump ($\vert [ \kappa ] (0) \vert \approx  0.1246$) in the curvature. However, this point  does not lie within the outflow boundary and hence it does not have any adverse effect on the rates of convergence. 
 There are  two external  characteristic points at $(0, \pm  M_y), M_y=2.5*\pi ^2+\beta $. Also $M_x = 2.25 \pi ^2+\beta $ and  the upper bound in (\ref{maxR}) is $R \approx < 14.28$.  
In Table \ref{Tab4} we present global orders   for  the  test problem, with $R=1$:  
\begin{equation}\label{Test3}
 -\ve \triangle u+ u_x +u=\Bigl( (1- \frac{y^2}{M_y^2})(\frac{x}{M_x})\Bigr)^4, \  (x,y) \in \Omega _2,\ 
u=0, \ (x,y) \in \partial \Omega _2 .
\end{equation}
 A  sample computed solution is displayed in Figure \ref{fig:ex2-solution}.

 \begin{table}[ht]
\centering\small
\begin{tabular}{|r|rrrrrrr|}\hline 
$\varepsilon| N$&$N=8$& 16&32&64&128&256&512\\ \hline 
1.000000 &   1.1815 &   0.8659 &   1.2460  &   0.7149  &   2.1452 &   0.4748 &   2.0765\\ 
0.500000 &   1.5030 &   0.9087 &   0.7154  &   1.0345  &   1.1353 &   1.4214 &   1.2746\\ 
0.250000 &   1.6797 &   1.0045 &   0.6723  &   1.1996  &   1.7662 &   1.0373 &   0.9828\\ 
0.125000 &   1.6394 &   1.2235 &   1.4583  &   0.9942  &   0.9491 &   0.9721 &   0.9782\\ 
0.062500 &   1.6173 &   1.3528 &   1.2847  &   0.9799  &   0.9524 &   0.9745 &   0.9955\\ 
0.031250 &   1.6061 &   1.3424 &   0.9507  &   0.9831  &   0.9934 &   0.9969 &   0.9986\\ 
0.015625 &   1.6006 &   0.4432 &   0.9738  &   0.9372  &   0.9678 &   0.9833 &   0.9913\\ 
0.007813 &   1.5978 &   0.6144 &   0.6812  &   0.8724  &   0.7878 &   0.7226 &   0.7270\\ 
0.003906 &   1.5964 &   0.6192 &   0.6815  &   0.8722  &   0.7880 &   0.7225 &   0.8377\\ 
0.001953 &   1.5957 &   0.6216 &   0.6816  &   0.8721  &   0.7881 &   0.7224 &   0.8377\\ 
0.000977 &   1.5954 &   0.6228 &   0.6816  &   0.8720  &   0.7881 &   0.7224 &   0.8377\\ 
0.000488 &   1.5952 &   0.6234 &   0.6817  &   0.8720  &   0.7881 &   0.7224 &   0.8377\\ 
0.000244 &   1.5951 &   0.6237 &   0.6817  &   0.8720  &   0.7881 &   0.7224 &   0.8377\\ 
0.000122 &   1.5951 &   0.6238 &   0.6817  &   0.8720  &   0.7882 &   0.7224 &   0.8377\\ 
0.000061 &   1.5951 &   0.6239 &   0.6817  &   0.8720  &   0.7881 &   0.7224 &   0.8377\\ 
0.000031 &   1.5951 &   0.6239 &   0.6817  &   0.8720  &   0.7881 &   0.7224 &   0.8377\\ 
0.000015 &   1.5951 &   0.6239 &   0.6817  &   0.8720  &   0.7881 &   0.7224 &   0.8377\\ 
0.000008 &   1.5951 &   0.6240 &   0.6817  &   0.8720  &   0.7881 &   0.7224 &   0.8377\\ 
0.000004 &   1.5951 &   0.6240 &   0.6817  &   0.8720  &   0.7881 &   0.7224 &   0.8377\\ 
0.000002 &   1.5951 &   0.6240 &   0.6817  &   0.8720  &   0.7881 &   0.7224 &   0.8377\\ 
0.000001 &   1.5951 &   0.6240 &   0.6817  &   0.8720  &   0.7881 &   0.7224 &   0.8377\\ \hline 
$\bar p ^N$ &   1.5951 &   0.4599 &   0.8440  &   0.8724 &   0.7878 &   0.7226  &   0.7270\\ \hline
\end{tabular}
\caption{Computed double-mesh global orders of convergence $\bar p_\ve ^N$ for the corrected approximations $\bar U$, using a strip of width $R=1$,
for problem (\ref{Test3}) and  $\beta =0.5$.}
\label{Tab4}
\end{table}

\noindent {\bf Example 3} The boundary is   $\partial \Omega _3:=\{(\phi(t), \psi (t)) \} = \{ (\rho \cos t, \rho \sin t) \}$, where $\rho := \beta + \cos ^2 t$ and  $0 < \beta <2$. 
The inflow boundary is disjointed and corresponds to the intervals
(in the $t$ variable)
\[
(0,\theta ) \cup (\pi/2,\pi-\theta ) \cup (\pi+\theta, 3\pi/2) \cup (2\pi-\theta, 2\pi), \quad  \theta  = \arcsin \sqrt{\frac{1+\beta}{3}}.
\]
 The domain  has four exterior characteristic points (where for $\beta =0.5, \kappa = 3/(2\sqrt{2})$) at \[
(\pm P, \pm Q), \quad P:= \frac{2(1+\beta)}{3}  \sqrt{\frac{2-\beta}{3}} \quad 
Q:=\frac{2(1+\beta)}{3}  \sqrt{\frac{1+\beta}{3}}\]
and two interior characteristic points (where for $\beta =0.5, \kappa = 2$) at $(0, \pm  \beta )$.
A test problem can be 
\begin{equation}\label{Test5}
 -\ve \triangle u+ u_x +u=(1-\frac{y}{\beta})^4(1+\frac{y}{\beta})^4H(y-\beta) H(y+\beta), \  (x,y) \in \Omega _3,\
u=0, \ (x,y) \in \partial \Omega _3;
\end{equation}
where $H(\cdot)$ is the Heaviside unit  step function.

 \begin{figure}[h!]
\centering
\resizebox{\linewidth}{!}{
	\begin{subfigure}[Domain $\Omega _3$ with $\beta =0.5$]{
		\includegraphics[width=7cm, scale =5]{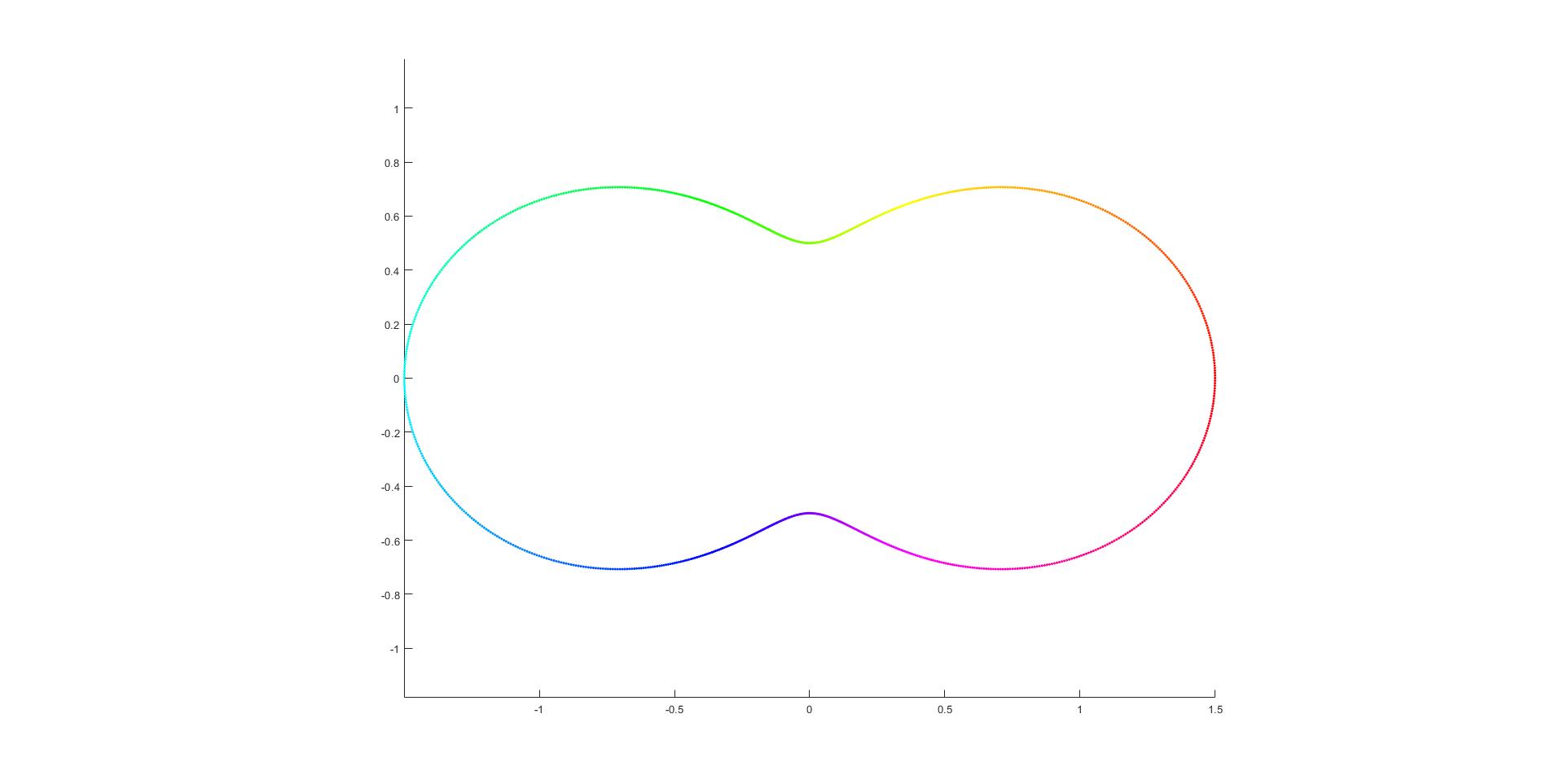} 
		}
    \end{subfigure}
\begin{subfigure}[Computed solution for $\ve =2^{-10}$]{
		\includegraphics[width=7cm, scale =5]{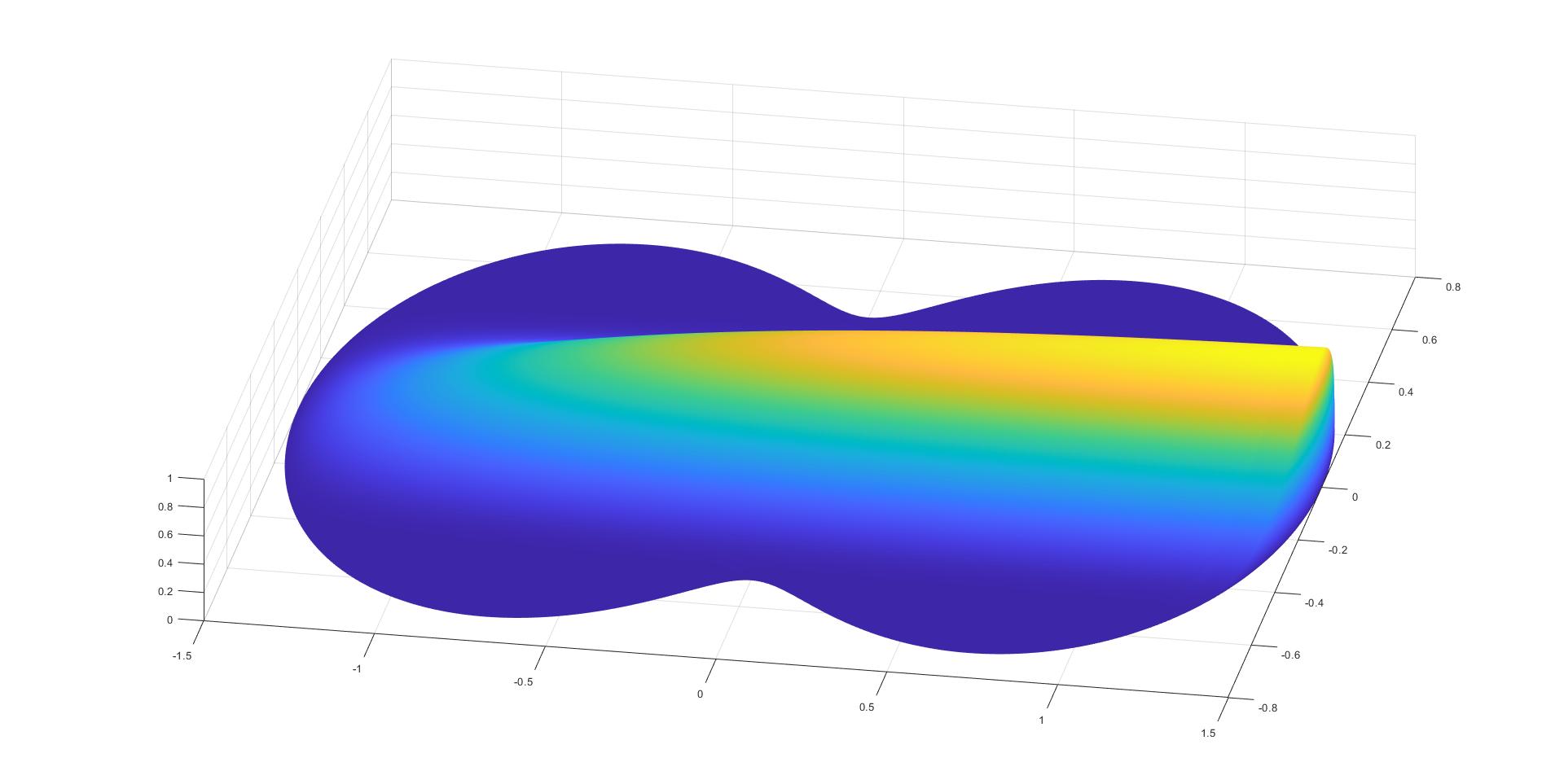} 
		}
	\end{subfigure}
}
	\caption{Computed solution of problem (\ref{Test5}) for $\ve =2^{-10}$ on the domain $\Omega _3$ with $\beta =0.5$}
	\label{fig:ex3-solution}
 \end{figure}

 \begin{table}[ht]
\centering\small
\begin{tabular}{|r|rrrrrrr|}\hline 
$\varepsilon| N$&$N=8$& 16&32&64&128&256&512\\ \hline 
1.000000 &   1.0702 &   0.3644 &   0.8105  &   0.5711  &   0.6534 &   0.8372 &   0.6988\\ 
0.500000 &   1.3472 &   0.6217 &   0.7877  &   0.5623  &   0.6173 &   0.8308 &   0.8403\\ 
0.250000 &   1.1339 &   1.5496 &   0.7490  &   0.4977  &   0.5958 &   0.8188 &   0.8354\\ 
0.125000 &   0.8425 &   1.5538 &   1.0672  &   1.1292  &   0.5590 &   0.7984 &   0.8272\\ 
0.062500 &   0.5542 &   1.5041 &   0.8783  &   1.1497  &   1.0228 &   1.2867 &   0.9587\\ 
0.031250 &   0.3465 &   1.8518 &   0.8686  &   1.0378  &   1.0317 &   1.2056 &   1.0415\\ 
0.015625 &   0.2300 &   2.2651 &   1.2117  &   1.4030  &   1.3229 &   0.7461 &   0.8992\\ 
0.007813 &   0.1682 &   2.6703 &   1.5264  &   1.3204  &   0.8332 &   0.7949 &   0.9250\\ 
0.003906 &   0.1367 &   2.8922 &   1.6560  &   1.0698  &   0.7922 &   0.8311 &   0.9404\\ 
0.001953 &   0.1206 &   2.9286 &   1.6488  &   1.0898  &   0.7699 &   0.8540 &   0.9416\\ 
0.000977 &   0.1124 &   2.8983 &   1.6037  &   1.1256  &   0.8283 &   0.8676 &   0.9377\\ 
0.000488 &   0.1082 &   2.8689 &   1.4319  &   1.2229  &   0.9477 &   0.8704 &   0.9388\\ 
0.000244 &   0.1003 &   2.8576 &   1.3268  &   1.2763  &   1.0112 &   0.8838 &   0.9416\\ 
0.000122 &   0.0930 &   2.8554 &   1.2726  &   1.2552  &   1.0361 &   0.9471 &   0.9417\\ 
0.000061 &   0.0892 &   2.8547 &   1.2456  &   1.2200  &   1.0735 &   0.9525 &   0.9553\\ 
0.000031 &   0.0873 &   2.8544 &   1.2323  &   1.2028  &   1.0944 &   0.9378 &   0.9365\\ 
0.000015 &   0.0863 &   2.8542 &   1.2257  &   1.1943  &   1.1063 &   0.9313 &   0.9236\\ 
0.000008 &   0.0859 &   2.8542 &   1.2224  &   1.1902  &   1.1128 &   0.9290 &   0.9170\\ 
0.000004 &   0.0856 &   2.8542 &   1.2208  &   1.1882  &   1.1161 &   0.9283 &   0.9142\\ 
0.000002 &   0.0855 &   2.8542 &   1.2200  &   1.1871  &   1.1179 &   0.9282 &   0.9132\\ 
0.000001 &   0.0854 &   2.8541 &   1.2196  &   1.1866  &   1.1187 &   0.9281 &   0.9127\\ \hline 
$\bar p ^N$ &   0.0854 &   2.8541 &   1.0716  &   1.0378 &   1.0317 &   1.2056  &   1.0191\\ \hline
\end{tabular}
\caption{Computed double-mesh global orders of convergence $\bar p_\ve ^N$ for the corrected approximations $\bar U$ 
for problem (\ref{Test5}) and  $\beta =0.5$.}
\label{Tab5}
\end{table}

\begin{remark}
The data in test problem (\ref{Test5}) has been chosen to satisfy Assumption 3, with $\delta =0.5$.
The data in the test problems (\ref{Test2}) and (\ref{Test3}) do not staisfy the compatibility constraints in Assumption 3 for any choice of $\delta >0$. Nevertheless, for all three test problems we observe parameter-uniform convergence in each of the corresponding Tables of orders of convergence. 
\end{remark}

\end{document}